\documentclass[12pt]{article}
\usepackage{amsmath,amsfonts,amssymb,amsthm}

\def\noi{\noindent}

\def\pf{\noi{\bf Proof.\ \,}}
\def\eop{{$\square$}}

\def\labtt#1{\label {#1}}
\def\labttr#1{\label {#1}\rm }

\def\o{\omega}

\def\vep{\varepsilon}

\def\DD{{\mathbb D}}
\def\FF{{\mathbb F}}
\def\QQ{{\mathbb Q}}
\def\RR{{\mathbb R}}

\def\ZZ{{\mathbb Z}}

\def\la{\langle}
\def\ra{\rangle}
\def\<{\langle}
\def\>{\rangle}

\def\half{{1 \over 2}}
\def\fourth{{1 \over 4}}

\def\dual#1{#1^*}        

\def\explus#1{2^{1+2#1}_+}

\def\ratholoex#1{2^{1+2#1}_+\Omega^+(2#1,2)}

\def\bw#1{BW_{{2^{#1}}}}

\def\brw#1{BRW^+(2^{#1})}

\begin{document}

\newtheorem{thm}{Theorem}[section]
\newtheorem{prop}[thm]{Proposition}
\newtheorem{lem}[thm]{Lemma}
\newtheorem{rem}[thm]{Remark}
\newtheorem{coro}[thm]{Corollary}
\newtheorem{conj}[thm]{Conjecture}
\newtheorem{de}[thm]{Definition}
\newtheorem{hyp}[thm]{Hypothesis}

\newtheorem{nota}[thm]{Notation}
\newtheorem{ex}[thm]{Example}
\newtheorem{proc}[thm]{Procedure}  
\newtheorem{app}[thm]{Application} 

\def\reftt#1{(\ref{#1})}

\def\threed{{\frak 3}}
\def\threedp{{\frak 3'}}
\def\oned{{\frak 1}}
\def\eightd{{\frak 8}}

\begin{center}

{\Large \bf Few-cosine spherical codes and Barnes-Wall lattices}  

\end{center}

\begin{center}
{\Large    }

Robert L.~Griess Jr.
\\[0pt]
Department of Mathematics\\[0pt] University of Michigan\\[0pt]
Ann Arbor, MI 48109-1109  \\[0pt] {\tt rlg@umich.edu}

\end{center}
 \centerline{28 April, 2006}
 \vfill \eject
\begin{abstract}Using Barnes-Wall lattices and 1-cocycles on 
finite groups of monomial matrices, 
 we give a procedure to
construct tricosine spherical codes.  
This was inspired by a 
14-dimensional 
code  which Ballinger, Cohn, Giansiracusa and Morris 
discovered in studies of the universally optimal property.  
It has 64 vectors and cosines $-\frac 3 7, -\frac 1 7, \frac 1 7$.  We construct the {\it Optimism Code}, a 4-cosine spherical code with 256 unit vectors in 16-dimensions.  The cosines are $0,\pm \fourth, -1$.  Its automorphism group has shape $2^{1+8}{\cdot}GL(4,2)$.  
The Optimism Code contains a subcode related to the BCGM code.  
The Optimism Code implies existence of a  nonlinear binary code with parameters $(16,256,6)$, a Nordstrom-Robinson code, and gives a context for determining its automorphism group, which has form $2^4{:}Alt_7$.
\end{abstract}
 
\tableofcontents 

\section{Introduction}

A {\it spherical code} is a finite set of unit vectors in Euclidean space.  A {\it cosine} of a spherical code is the inner product of distinct unit vectors in the code.

\begin{de}\labttr{cossc}  
Call a spherical code {\it $n$-cosine} if the inner products of distinct unit vectors form an $n$-set.  
When $n=3$, we use the term  {\it tricosine}. 
\end{de}

We present a general existence criterion \reftt{unidefsc} for tricosine spherical codes, based on the unidefect concept.  We record an infinite series of examples and some special ones in dimensions 14 to 16.   The 14-dimensional one probably satisfies the hypotheses of 
 the recent uniqueness result  \cite{bb} and if so would be isometric to  
 $\mathcal {BCGM}$.  

We construct and analyze the {\it Optimism Code}, a spherical 4-cosine $(16,256,6)$ code which can be used to derive all our special examples.  Its existence depends on an easy result from group extension theory.  The Optimism Code has isometry group a nonsplit extension $2^{1+8}GL(4,2)$ and is connected to a nonlinear binary Nordstrom-Robinson type code.  For the latter code, we have an easy existence proof and determination of its automorphism group.  
We have two existence proofs for both the 14-dimensional code and the Optimism Code, both involving 1-cocycles (also called derivations; see \reftt{rightderiv}).

This article was inspired by 
a spherical code $\mathcal {BCGM}$  found by 
Brandon Ballinger, Henry Cohn, Noah Giansiracusa and Elizabeth Morris,  while investigating the universally optimal property \cite{ck}.

{\bf BCGM1.}  $\mathcal {BCGM}$ has 64 unit vectors in dimension 14, two of which make angles with cosines $\{-\frac 37, -\frac 17, \frac 17 \}$.  

{\bf BCGM2.}  Its isometry group $\widetilde H$ has these properties: 

(i)  $O_2(\widetilde H)$ is nonabelian of order $2^7$, $O_2(\widetilde H)'=Z(O_2(\widetilde H))$.

(ii) $\widetilde H /O_2(\widetilde H)\cong GL(3,2)$.

{\bf BCGM3.} $\mathcal {BCGM}$ is an association scheme.  

\bigskip 

The rhythm of  
$\{-\frac 37, -\frac 17, \frac 17 \}$ 
suggested the ZOPT property for Barnes-Wall lattices \cite{bwy}.    These lattices were therefore considered a possible source for interesting spherical codes. 
For $\mathcal {BCGM}$, it is natural to look at $\bw 4$, whose automorphism group $\brw 4\cong 2^{1+8}_+\Omega^+(8,2)$ contains many subgroups which look roughly like $\widetilde H$.

Notation and terminology follows that in \cite{bwy} and \cite{ibw1}.  
The {\it cubi theory}  of \cite{ibw1} is recommended (see Section 3, especially 3.21 ff.).  
A few techniques from group cohomology are collected in an appendix.

\begin{table}[htdp]
\caption{Brief summary of our unidefect tricosine codes }
\begin{center}
(In row 1, $3\le m \le d$ and $2^m-1$ is a  Mersenne prime.) 
\vskip 0.5 cm 
\begin{tabular}{|c|c|c|c|}
\hline
Symbol & Dimension & Number of unit vectors & Cosines \cr
\hline
\hline
$\mathcal {DSC}_{2^d-\ell, 2^{d+m}}$ & 
$2^d-\ell$, $\ell$ small & $2^{m+d}$ & $\frac {-2^{d-k}-\ell}{2^{d}-\ell}, \frac {-\ell}{2^{d}-\ell}, \frac {2^{d-k}-\ell}{2^{d}-\ell}$\cr 
\hline
$\mathcal {NSC}_{16, 64}$ & 
16 & 64 & $-\frac 1 4 , 0 , \frac 1 4$ \cr 
\hline 
$\mathcal {NSC}_{15, 64}$ & 
15 & 64 & $-\frac 15, -\frac 1{15}, \frac 13$ \cr 
\hline 
$\mathcal {NSC}_{14, 64}$ & 
14 & 64 & $-\frac 37, -\frac 17, \frac 17$ \cr 
\hline 
$\mathcal {NSC}_{16, 128}$ &
16 & 128 & $-\frac 1 4 , 0 , \frac 1 4$ \cr 
\hline 
$\mathcal {NSC}_{15, 128}$ & 
15 & 128 & $-\frac 15, -\frac 1{15}, \frac 13$ \cr 
\hline 

\hline 
\end{tabular}
\end{center}
\label{default}
\end{table}

{\bf Acknowledgements.  }  We thank Henry Cohn for describing $\mathcal {BCGM}$ and explaining background.  For useful consultations, we thank 
Eichii Bannai, Etsuko Bannai and Akihiro Munemasa.  
This work was begun at the Oberwolfach Mathematische Forschungsinstitut at the meeting 21-25 November, 2005, and was supported in part by NSA grant USDOD-MDA904-03-1-0098.

\begin{nota}\labttr{general}
{}\ 

$A.B$ means an extension of groups ($A$ normal, giving quotient $B$).  

$A{\cdot}B, A{:}B$ mean nonsplit, split extensions, respectively.
  
$p^m$ means an elementary abelian $p$-group of rank $m$ ($p$ prime).

$O_p(G)$ means the largest normal $p$-subgroup of the finite group $G$ ($p$ prime).

$C_G(X), N_G(X)$ shall mean the centralizer, normalizer, respectively, in a group $G$ of a subset $X$ (subscript $G$ may be omitted); this notation extends to subsets $X$ of a set on which $G$ has a permutation representation.

$P(S)$ means the power set of the set $S$ and $PE(S)$ means the subspace of even subsets of the finite set $S$.  

The term {\it weight} refers to weight of a binary codeword, i.e., the cardinality of its support.  We generally identify a binary codeword with its support and vector addition with the symmetric difference of subsets.  

$Mon(n, \{\pm 1\})$ denotes the group of degree $n$ monomial matrices with entries $0, \pm 1$ only.

  Groups actions on sets and modules will be on the right, sometimes with exponential notation.  The conjugate of $x$ by $y$ is $x^y=y^{-1}xy$ and the commutator of $x$ and $y$ is $[x,y]=x^{-1}y^{-1}xy$.  
  
    \end{nota}

\section{Unidefect criterion for a tricosine spherical code}

We give a criterion for constructing tricosine spherical codes in \reftt{unidefsc}.  First, we need to sketch notations for the Barnes-Wall lattices, $\bw d$.  This is taken from our recent article \cite{bwy}.  See also the classic articles \cite{bw}, \cite{be}.

\begin{nota} \labttr{bwsetup} 
We take the rank $2^d$ Barnes-Wall lattice $\bw d$ and the subgroup $G:=\brw d \cong \ratholoex d$ of the automorphism group, which for $d\ne 3$ is the full automorphism group.  
As usual, $R=O_2(G)\cong \explus d$.   
Let $F$ be a sultry frame and $\mathcal B \subset F$, an orthogonal basis of $V$, the ambient real vector space.  
We use indices $\{0,1,\dots, 2^d-1\}$ to  
label the orthonormal basis $2^{- \frac d 2}\mathcal B= \{ v_0, v_1, \dots \}$ 
and vector space $\Omega:=\FF_2^d = \{\o_0, \o_1, \dots \}$.   
When $A$ is a subset of $\Omega$, write $v_A:=\sum_{i \in A} v_i$.

For a subset $A$ of $\Omega$, let $\vep_A$ be the orthogonal transformation which takes each
$v_i$ to $\begin{cases} -v_i & i \in A\cr v_i & i \not \in A \end{cases}$.   
 In $G$, take the associated diagonal subgroup $D$ and $N$ its normalizer.  The group $D$ consists of all $\vep_A$ where $A$ ranges over the Reed-Muller code $RM(2,d)$.  Also  we assume that $\mathcal B$ is chosen so that 
 $N$ is a semidirect product $DP$, for a group $P\cong AGL(d,2)$ of permutation matrices with respect to $\mathcal B$.   We assume that the bijection $v_i \mapsto \o_i$ is an equivariance respecting the identification 
$P\cong AGL(d,2)$.  
\end{nota}

\begin{nota}\labttr{trivscsetup}  We identify $V$ with $\RR^{2^d}$ by use of  the orthonormal basis $v_i, i \in \Omega$ \reftt{bwsetup}.  
Let $x_0=(1,1,1, \dots, 1,1)=v_\Omega$ be the all-$1$ vector with respect to the orthonormal basis $v_i, i \in \Omega$; if 
$d$ is even, it is in the standard $\bw d$. 
Suppose that $Q$ is a subgroup
of $P$ and $J$ a subgroup of $D$ which is normalized by $Q$ and $-1 \not \in J$.  
\end{nota}

\begin{de}\labttr{orbitcode}
A spherical code is a {\it diagonal code} if it is an orbit of $x_0$ by a subgroup of the diagonal group $D$ in a $\brw d$-group.  
\end{de}

\begin{rem}\labttr{aboutdefect} 
The {\it defect} of an involution $t\in G$ is the integer $k$ so that 
$2^{2(d-k)}$ is the order of $C_{R/Z(R)}(t)$.  At once, $0\le k \le \frac d 2$.  
See \cite{ibw1} for discussion of defect of an involution in $G$ and of a codeword in $RM(2,d)$.  The main properties we need here are that the weight of a defect $k$ codeword 
is one of  the values 
$2^{d-1}$ or 
$2^{d-1}\pm 2^{d-k-1}$,  and the fact that every involution of $G$ in a given coset of $R$ has a common defect (this implies that  defect is constant for any coset of $RM(1,d)$ in $RM(2,d)$).    A codeword of weight $2^{d-1}\pm 2^{d-k-1}$ is {\it clean} and one of weight $2^{d-1}$ is {\it dirty}.  
\end{rem} 

\begin{de} \labttr{unidefcond} 
Call a function $f : Q \rightarrow J$
a {\it near-derivation} if the associated function
$\bar f : Q \rightarrow J/[J\cap R]$ is a derivation, i.e., a 1-cocycle; see \reftt{rightderiv}.  
The {\it strong unidefect condition} on the function $f : Q \rightarrow J$ 
is that 
there exists a fixed integer 
$k$, $1\le k \le \frac d 2$  so that 
for all $x\in Q$, $f(x)\in R$ or $f(x)$ 
is a defect $k$ involution, i.e., has 
the form $\vep_A$, where
$A$ is a defect $k$ codeword. 
An alternate formulation  is that 
there exists a fixed integer 
$k$, $1\le k \le \frac d 2$  so that 
every value of $f$ is in $R$ or 
has defect $k$.  

The {\it unidefect condition} on $f$ is that every nonidentity value of 
$f$ is an involution of trace 0 or of 
the form $\vep_A$, where
$A$ is a clean defect $k$ codeword.   This condition is weaker than strong unidefect because it allows $f(x)$ to have defect not $k$ as long as $f(x)$ has trace 0.  
\end{de}

\begin{prop} \labtt{unidefsc}  
Let $Q, J, f$ be as in \reftt{unidefcond}.  
Assume that $f$ satisfies the  unidefect condition for defect $k$.  Let $H:=(J\cap R)\{f(x)x|x\in Q\}$ be the group containing $J\cap R$ which is associated to $f$; see \reftt{interp1coc}.   

Then the orbit $x_0H$ is a set of vectors of common norm $2^d$ for which the inner product of two distinct members is 0 or $\pm 2^{d-k}$.   
The length of this orbit is $|J\cap R||Q{:}Ker(\bar f)|$ (see \reftt{kerder}).  
\end{prop}

\begin{rem}\labttr{aboutunidef} 
(i) The ZOPT property of Barnes-Wall lattices is that two minimal vectors have inner product which is 0 or plus or minus a power of 2.

(ii) When $f=0$, the code is diagonal \reftt{orbitcode}.   

(iii) 
A change in defect of the cocycle may cause a  change in the cosines.  
\end{rem} 

\begin{de} \labttr{reduction}  
Suppose that $S$ is a set of equal norm nonzero vectors in $V$ and that $W$ is a subspace of $V$ so that every element of 
$S$ has the same projection to $W^\perp$.  Then $S$ may be projected to $W$ and rescaled to make a spherical code in $W$.  If $S$ is $n$-cosine, then so is the projection.   This process is called {\it reduction to $W$} and the resulting code is called the {\it reduction of $S$} or just the {\it reduced code}.
\end{de} 

\begin{prop}\labtt{reducedcosines}
In the situation of \reftt{unidefsc}, \reftt{reduction}, denote $\ell :=  dim(W^\perp)$.  The reduced spherical code has cosines 
$\frac {-2^{d-k}-\ell}{2^{d}-\ell}, \frac {-\ell}{2^{d}-\ell}, \frac {2^{d-k}-\ell}{2^{d}-\ell}$.  
\end{prop} 
\pf
The projection of the norm $2^d$ vectors of $x_0H$ to $W^\perp$ have norm $2^d-\ell$. 
\eop

\begin{de} \labttr{unidefectcode}  
A spherical code is called a {\it unidefect code} if it is 
created from an orbit as in \reftt{unidefsc} by scaling to unit length; or created from an orbit by projecting and rescaling as in \reftt{reduction}.  
\end{de}  

\section{Diagonal codes}  

\subsection{$\mathcal {DSC}_{2^d-\ell, 2^{d+m}}$, for small $\ell$}

Even the case of small $Q$  \reftt{trivscsetup} and $f=0$ \reftt{aboutunidef}  
is interesting, as the following examples show.

\begin{de}\labttr{pure}  Fix an integer $k>0, k\le \frac d 2$.  
A subset $Y$ of $D$ is {\it defect $\{0,k\}$-pure} if every involution in it has defect 0 (i.e., is in the lower group $R$) 
or defect $k$.  
\end{de}

It would be useful to find large pure {\it subgroups}.

\begin{de}\labttr{mersenneex}  {\it A family of diagonal codes associated to Mersenne primes.  }

Let $p=2^m-1$ be a Mersenne prime 
and suppose that $3\le m \le d$.  Take $g\in P$ of order $p$ and assume that $g$ fixes $v_0$.  For the action of $Q:=\la g \ra$ on $D$, every irreducible constitutent has dimension 1 or $m$.
When $m=d$, there is a single nontrivial constituent in $D\cap R$ and $\frac {d-1}2$ of them in $D/D\cap R$.

Let $J$ be a $\la g \ra$-invariant subgroup of $D$ so that 
$J$ fixes $v_0$, 
$J\cap R=C_{D\cap R}(v_0)$ and 
$J/J\cap R$ is a $\la g \ra$-irreducible module.  
Then $J\cap R$ has order $2^d$ and $J$ has order $2^{m+d}$.  
Since $\la g\ra$ acts transitively on the nontrivial elements of 
$J/J\cap R$, 
$J$ is $\{0, k\}$-pure for some $k>0$.  

Now take the orbit $x_0J=x_0\la g\ra J$, which is in bijection with $J$ so has $2^{m+d}$ elements.  The inner product of $x_0$ with any other member of this orbit is one of $0, \pm 2^{d-k}$.  Transitivity implies the analogous property for every member of the orbit.  Rescaling gives unit vectors with inner products $0, \pm 2^{-k}$.   
\end{de}  

Every nonidentity element of $J\cap R$ has trace 0.  If $r$ is such an element $(x_0, x_0r)=0$.  
For completeness, we record the following.  

\begin{lem}\labtt{alloccur}
Both $\pm 2^{d-k}$ 
occur as inner products in the situation of \reftt{mersenneex}.  
\end{lem}  
\pf 
We use the orthogonality relations for characters of $J$.  
Let $\chi$ denote the trace function for linear transformations on $V$ \reftt{bwsetup}.

First we show that a nonzero inner product occurs.  Assume otherwise.  
Then $s:=\sum_{y \in J} \chi (y)= \sum_{y \in J\cap R} \chi (y)$, which is $2^d$ since $V|_{J\cap R}$ is the regular representation.  This is a contradiction since $s$ is divisible by $|J|=2^{d+m}$.  

Second, $\frac s {2^{d+m}}$ is the multiplicity of 1 in $V|_J$, which is at least 1 since $J$ fixes $v_0$ and is at most 1 since $V|_{J\cap R}$ affords the regular representation of $J \cap R$.  So, $\frac s {2^{d+m}}=1$, whence $s=2^{d+m}> 2^d$.  It follows that the unique positive candidate $2^{d-k}$ occurs as an inner product.  

Thirdly, we must show that  $-2^{d-k}$ occurs.  Suppose not.  Then $\chi (y)\ge 0$ for all $y \in J$.  
Let $h \in J \setminus R$.   In $(J\cap R)h$, the number of clean elements is $2^{2k}$ \cite{ibw1}, Prop. 3.32.  Therefore 
$2^{d+m}=s=\sum_{y \in J} \chi (y) =2^d+2^{2k}2^{d-k}(2^m-1)=2^d+2^{d+k}(2^m-1)$.  Dividing both sides by $2^d$ gives $2^m=1+2^k(2^m-1)$, which is impossible since the right side is odd.  
\eop

\begin{nota}\labttr{restricteddsc}
The spherical code of \reftt{mersenneex} is denoted 
$\mathcal {DSC}_{2^d, 2^{m+d}}=\mathcal {DSC}_{2^d, 2^{m+d};J,g}$.  The notation 
$\mathcal {DSC}_{2^d-\ell, 2^{m+d}}=\mathcal {DSC}_{2^d-\ell, 2^{m+d};J,g,W}$.
 means a 
spherical code obtained by 
projecting $\mathcal {DSC}_{2^d, 2^{m+d};J,g}$ to $W$, the orthogonal of an $\ell$-dimensional space fixed pointwise by $J\la g \ra$, e.g., $span\{v_0\}$.  
\end{nota}

\begin{ex}\labttr{diagex} 
For example, take $d=m=5, p=31$.  We get a tricosine spherical code of  1024 elements in $\QQ^{32}$ in which nonzero inner products are either  $\pm \half$ or $\pm \fourth$.  One may project to a 31-space for another code with respective cosine set either $\{-\frac {17} {31}, -\frac 1 {31}, \frac {15}{31}\}$ or $\{-\frac 9{31}, -\frac 1{31}, \frac 7{31}\}$.  
\end{ex}

\begin{rem}\labttr{whichk} 
(i) 
It is not obvious which $k$  in the range
$1\le k \le \lfloor \frac d 2 \rfloor$
may occur  this way.

(ii) 
 The full isometry group of such a code could contain $J\la g \ra$ properly, e.g. \reftt{sphcode1} occurs here for $p=7, m=d=3$.  
 
 (iii) 
 When $p=7, d=4$, there are several (possibly nonisometric) codes, all with defect 1.  
\end{rem}

\section{Nondiagonal codes}

We construct nondiagonal spherical codes in dimension 14 through 16 using 1-cocycle theory for  the simple group of order 168 acting on various sections of the frame stabilizer.

\begin{nota}\labttr{14setup} 
We continue to use the notation of \reftt{bwsetup} for $d=4$. 
Let $W$ be the orthogonal of $W_{01}:=span\{v_0, v_1\}$.
We take the subgroup $E\cong 2^3$ of $D\cap R$ which is trivial on $W_{01}$.  
Let $P_{\{01\}}\cong 2 \times 2^3{:}GL(3,2)$ be the subgroup of $P\cong AGL(4,2)$ which stabilizes $\Omega_{\{01\}}=\{\o_0,\o_1\}$ and 
$P_{0,1}$ the subgroup which fixes both $\o_0$ and $\o_1$.  
It has the form $P_{0,1}=UQ$, where $Q \cong GL(3,2)$ and $U \cong 2^3$.  
We take $B_0$ to be an affine hyperplane of $\Omega$ which contains $\o_0$ but not $\o_1$ 
and let $B_1$ be its complement,  
an affine hyperplane of $\Omega$ which contains $\o_1$ but not $\o_0$.  We choose $Q$ to stabilize $B_0$ and $B_1$.  
Let $p_{01}$ be the involution which generates $Z(P_{\{01\}})$.  It corresponds to translation from $\o_0$ to $\o_1$.  

We use the Reed-Muller codes $RM(r,4)$ spanned by affine subspaces of $\Omega$ of codimension $r$.  
Let $S_i$ be the subspace of $P(\Omega)$ which is spanned by all affine codimension 2 subspaces which are contained in $B_i$.  Then $dim(S_i)=4$ for $i=0,1$, $S_0\cap S_1=0$, 
$RM(2,4) \ge S:=S_0+S_1\ge RM(1,4) \ge S_{01}:=C_{S}(p_{01})$ and 
$dim(S_0+S_1)=8$.

Define $T_i:=\{A \in S_i|\o_i \not \in A\}$, a dimension 3 subspace of $S_i$, for $i=0,1$.  Let $F_i:=\{\vep_A|A\in T_i\}\cong 2^3$.  
Define $T_{01}:=C_{T_0+T_1}(p_{01})=\{ x+x^{p_{01}} |x \in T_0\}$.

For $\{i,j\}=\{0,1\}$, define $D_i:=\la \vep_A | A \in PE(B_i), \o_i \not \in A \ra \cong 2^6$.  The groups $D_i$ are not contained in $Aut(L)$ (in fact, $D_i\cap Aut(L)=F_i, i=0, 1$) but  the diagonal group  $D_{01}:=C_{D_0\times D_1}(p_{01})$ is in $Aut(L)$.  
The corresponding subspace of $PE(B_0)\oplus PE(B_1)$ is denoted  
$B[01]$, so that $D_{01}=\{\vep_A|A\in B[01]\}$.  
Moreover,  
$J:=F_0F_1D_{01}=C_D(\{v_0, v_1\})$ is an index 4 subgroup of $D$.  Also, $E=D_{01}\cap J$.

The action of $Q$ fixes the $S_i$ and $T_i$ and $Q$ normalizes the $F_i$ and $D_{01}$, so that $J=F_i \times D_{01}$ as a $Q$-module, for $i=0, 1$.  
\end{nota}

\subsection{$\FF_2[GL(3,2)]$-modules}

\begin{nota}\labttr{irrmods} 
Call the irreducible $\FF_2[GL(3,2)]$-modules $\threed, 
\threedp$, $\frak 1$ and $\frak 8$ (the number indicates dimension and the prime indicates duality). We inflate this notation to $\FF_2Q$-representations.  
Let us say that $\Omega \cong \FF_2^4$ as a $Q$-module has composition factors $\frak 1, \threed$.  
\end{nota} 

\begin{lem}\labttr{whichirreds} 
(i) 
$U\cong E \cong \threedp$.  

(ii) 
$T_0 \cong T_1 \cong \threedp$.  

(iii) $F_0\cong F_1\cong \threedp$.

(iv) 
$D_0/F_0 \cong D_1/F_1 \cong J/F_0F_1 \cong \threed$.  
\end{lem} 
\pf
(i) If we take $\o_0$ as an origin, 
$U$ is in $\FF_2$-duality with the quotient space 
$\Omega/\{\o_0,\o_1\}$.

(ii) The first isomorphism is realized by the action of $p_{01}$.  For the second, note that $T_1$ may be identified with linear functionals on $\Omega$ which have $\{\o_0,\o_1\}$ in their kernel.  

(iii) Consider the definition of $F_i$.

(iv) 
First, note that each $D_i/F_i$ is in duality with $T_i$.   Secondly, note that 
each $D_i/F_i$ covers $J/F_0F_1$.   
\eop

\subsection{Good subgroups of shape $2^3{.}GL(3,2)$}

\begin{nota}\labttr{good168}  
We continue to use the notations of \reftt{14setup}.  
We form the semidirect product $JQ$ and consider 
$\mathcal G:= \{ H | E \le H \le JQ, H/E\cong GL(3,2) \}$.  
This set 
has cardinality 256 and 
is a union of four orbits under $J$ or $JQ$, 
by \reftt{whichirreds}, \reftt{1coh168}.  
Each orbit is represented by a 1-cohomology class of $Q$ with coefficients in $J/E$.  
If $\gamma$ is a near-derivation \reftt{unidefcond} associated to $H$, write $\gamma=\gamma_0\gamma_1$ to indicate the components with respect to the direct sum $F_i \times D_{01}$, for a fixed $i\in \{0,1\}$.   Then $\gamma_0$ is a derivation and $\gamma_1$ is a near-derivation.  
Write $\bar \gamma_j$ for values 
of $\gamma_j$ 
 modulo $E$ and $\bar \gamma$ for values of $\gamma$ modulo $E$.    
Because of the correspondence of $H \in \mathcal G$ with the class of a near-derivation on $Q$, we may say that  $H$ has the unidefect property \reftt{unidefcond} if and only if such a near-derivation does.  
\end{nota} 

\begin{lem}\labtt{interpgam0gam1} 
Assume the notation as in \reftt{good168}.  
Then $H$ splits over $E$ if and only if $\bar \gamma_1$  is an inner derivation.  
\end{lem} 
\pf  The ``only if" part is trivial.  Assume $\bar \gamma_1$  is an inner derivation.  Then $EH$ is conjugate by an element of $D_{01}$ to $F_iEQ$ modulo $F_i$, which is split over $E$.  
\eop

\begin{rem}\labttr{orbitlengths} 
For each $H\in \mathcal G$, the orbit $x_0H$ is a spherical code whose set of cosines depends on  $H$.  We get its cardinality from \reftt{unidefsc} and the observation that the stabilizer in $H$ of $x_0$ is just $H \cap P$.  In the notation of \reftt{good168}, $H\cap P=H\cap Q= Ker(\bar \gamma_0)\cap Ker(\bar \gamma_1 )$.  The derivation kernels can have indices 1, 7 or 8 in $Q$ and $Ker(\bar \gamma_0)\cap Ker(\bar \gamma_1 )$ can have indices 1, 7, 8, 42 or 56.

Note that the orbit lengths depend on actual cocycles and not just cohomology classes.  
We are looking for a code like $\mathcal {BCGM}$, so the case of interest is $|Q:Ker(\bar \gamma_0)\cap Ker(\bar \gamma_1 )|=8$, which means 
$Ker(\bar \gamma_0)=Ker(\bar \gamma_1 )$ is a Frobenius group of order 21.  Derivations on irreducible 3-dimensional modules with such  kernels are outer and furthermore are associated to nonsplit extensions \reftt{interpgam0gam1}.  

Finally, we comment that the orbit  corresponding to a split extension contains groups of defects 0 and 1 only.  
\end{rem}

\begin{lem}\labtt{if2only2} 
Suppose that $\bar \gamma_0$ and $\bar \gamma_1$ have the same kernel and $\bar \gamma$ takes a nontrivial value which has defect $k\ge 1$.  Then all nontrivial values of $\bar \gamma$ have defect $k$ 
and so $Q, \gamma, J$ satisfy the unidefect $k$ condition.  
\end{lem} 
\pf 
We have $K:=Ker(\bar \gamma )=Ker(\bar \gamma_0)=Ker(\bar \gamma_1)$ is isomorphic to $Sym_4$ or a Frobenius group of order 21.  The action of $Q$ on the cosets of $K$ is doubly transitive.  Now use \reftt{propderiv}(iv).  
\eop

\begin{lem}\labtt{bothdefects}
Both unidefect 1 and 2 subgroups occur in $\mathcal G$.  
In particular, the class with both components \reftt{good168} noninner has a unidefect 2 subgroup $H^*$  which furthermore has orbit length $|x_0H^*|=64$.  
\end{lem}

\pf
A near-derivation of $Q$ with coefficients in one of the $F_i$  is a derivation.   If nontrivial, the derivation takes values which 
are involutions of defect 1 (because all nonzero codewords of $T_i$ have weight 4).  
A member of $\mathcal G$ with  
$\bar \gamma_1$ trivial has  defect 1 or 0.  

Consider the case $\bar \gamma_1$ nontrivial and $\bar \gamma_0$ trivial.  
Every weight in $B[01]$ is divisible by 4 and 
all codewords in $B[01]\setminus RM(1,4)$ have defect 1.  
Therefore, $\gamma$ has unidefect 1.

Suppose  $\bar \gamma_1$ noninner and $\bar \gamma_0$ noninner.  
Assume further that $Ker(\bar \gamma_0)=Ker(\bar \gamma_1)$, whence both are Frobenius groups of order 21 \reftt{anatomy}(ii).  This equality does occur for some groups in this orbit.  
We shall demonstrate explicitly such a  
$\gamma$ which 
takes value in $E\vep_Y$, for a 6-set $Y$.

Fix an involution $t\in Q$. 
 
Note that on $B_0\setminus \{v_0\}$, the action of the involution $t\in Q$ has a pair of length 2 orbits, hence  on $PE(B_0\setminus \{v_0\})$ has 2-dimensional commutator space, $M$.    Let $\{a,a'\}$ and $\{b,b'\}$ be the nontrivial orbits of $t$.  Then $\{a,a'\}, \{b,b'\}$ span $M$ and the 1-space $M \cap T_0$ is the span of $\{a,a',b,b'\}$, an affine 2-space. 

There are fixed points $d, e$ of $t$ on $B_0\setminus \{v_0\}$ 
so that $\{a,a',d,e\}$ is an affine 2-space (so is in $T_0$ \reftt{14setup}). 

Similarly,  $b, b'$ is contained in $\{b,b',d,e\}$, an affine 2-space which is the sum of the two previous affine 2-spaces.

Let $f$ be the remaining fixed point.  The 4-set 
$\{a,a',d,f\}$ is congruent to $\{e,f\}$ modulo $span(\{a,a',d,e\})$.  Both these sets are fixed by $t$.  
 
Now define $u:=e^{p_{01}}, v:=f^{p_{01}}\in B_1$.
The 6-set $Y:=\{u,v,a,a',d,f\}$ is fixed by $t$ and  
$Y = \{u, v, e, f\}+\{a, a', d, e\} \in B[01] + T_0 \le RM(2,4)$.  
There exists a near-derivation $\gamma$ so that $\gamma (t)=\vep_Y$ (see \reftt{involext}(ii), applied to $B[01]/T_{01}$ and $T_i$).  Now use \reftt{if2only2}.  
\eop

\begin{coro}\labtt{summary1} 
There are $H\in \mathcal G$ which have unidefect 2.  For such $H$,  
let $\gamma$ be an associated near-derivation.  
Every element of $H$ has the form $rq$, where $q$ is a permutation matrix and $r\in D$ effects sign changes at no coordinates, or at a clean codeword of defect $2$ 
(of weight 6 or 10) 
or at a midset (weight 8).   The values of $\gamma$ outside $R$ 
have defect 2.  
The extension $1\rightarrow E\rightarrow H \rightarrow GL(3,2)\rightarrow 1$ does not split.  There exists a particular such $H$, called $H^*$, so that $p_{01}$ normalizes $H^*$ and satisfies $[H^*,p_{01}]=O_2(H^*)$
\end{coro}  
\pf  All is clear except possibly for the nonsplitting.  For that point, we use \reftt{interpgam0gam1}, \reftt{anatomy}(ii) and the fact that $Ker(\bar \gamma_0)=Ker(\bar \gamma_1)$ is a Frobenius group of order 21.  
\eop

\subsection{$\mathcal {NSC}_{16, 64}$, $\mathcal {NSC}_{15, 64}$ and $\mathcal {NSC}_{14, 64}$}

\begin{nota}\labttr{x0}   
Let $H^* \in \mathcal G$ be a unidefect 2 group, as in \reftt{summary1}.  Then $x_0H^*=x_0\la p_{01}, H^* \ra$ has cardinality 64.  

  Let $\pi$ be the orthogonal projection $V\rightarrow W$ and $\rho$ the orthogonal projection $V\rightarrow  v_0^\perp$.  
 We define spherical codes 
$\mathcal {NSC}_{16, 64}$, $\mathcal {NSC}_{15, 64}$, $\mathcal {NSC}_{14, 64}$ as the vectors of the respective orbits 
$x_0H^* $, $(x_0H^*)\rho$ and $(x_0H)\pi $, scaled to be unit vectors in 16-, 15- and 14-dimensional space.  
\end{nota}

\begin{thm}\labtt{sphcode1}  (i) 
The cosines for $\mathcal {NSC}_{16, 64}$ are just  $0, \pm \fourth $.  

(ii) 
The cosines for $\mathcal {NSC}_{15, 64}$ are just $\{-\frac {1}5, -\frac 1{15}, \frac 13\}$.  

(iii) 
The cosines for $\mathcal {NSC}_{14, 64}$ are just $\{-\frac {1}7, -\frac {3}7, \frac 1 7\}$.  
\end{thm}
\pf \reftt{reducedcosines}. \eop

\subsection{$\mathcal {NSC}_{16,128}$ and $\mathcal {NSC}_{15,128}$}

Recall the definitions of $B_0$ and $B_1$ \reftt{14setup}.  
Consider, in $\bw 4$, the sign change isometry $\vep_{B_0}$, where $\o_0\not \in {B_0}$ and $\o_1\in {B_0}$.  

\begin{nota}\labttr{hstar} 
We increase $H^*$ to $H^{**}:=H^*\la \vep_{B_0}\ra$.  
Since $[H^*, \vep_{B_0}]\le E$, the dihedral group $\la p_{01},\vep_{B_0}\ra $ normalizes $H^*$ but $\vep_{B_0}$
does not normalize the group $H^*\la p_{01} \ra$ of \reftt{summary1}.
Denote by 
$\mathcal {NSC}_{16,128}$ the spherical code in $\RR^{16}$ obtained by scaling  the elements of $x_0H^{**}$ to unit length.  
Denote by 
$\mathcal {NSC}_{15,128}$ the spherical code in $\RR^{15}$ obtained projecting $x_0H^{**}$ to $v_0^\perp$, then 
rescaling to unit length.   
\end{nota}

\begin{thm}\labtt{128codes} 
(i) 
$\mathcal {NSC}_{16,128}$ has cardinality 128 and cosines $\{-\fourth, 0, \fourth \}$; 

(ii) $\mathcal {NSC}_{15,128}$ has cardinality 128 and cosines $\{-\frac {1}5, -\frac 1{15}, \frac 13\}$. 
\end{thm} 
\pf  As with $H^*$, use 
\reftt{unidefsc} and 
\reftt{reducedcosines}.  Since $H^*$ satisfies  the strict 2-unidefect condition, so does 
$H^{**}$, which is created from $H^*$ by  
replacing $E$ with 
the slightly larger lower group $E\la \vep_{B_0} \ra$. 
\eop

\begin{rem}\labttr{rem128}   
(i) The automorphism group of $x_0H^{**}$ excludes $p_{01}$ (or else $-1=[p_{01},\vep_{B_0}]$ would be an automorphism, which does not respect the cosine set).

(ii) 
Projection to the 14-space $W$ does not seem to give a tricosine code.  
\end{rem}

\section{Computations}  

We outline a straightforward computational method for finding our  unidefect spherical codes \reftt{unidefsc} by computer.   
Such a code is an orbit for a finite group $H$, so is a union of orbits of any subgroup of $H$.  We take the subgroup $H\cap P$, the subgroup of $H$ consisting of permutation matrices.  Since there are no signs in these matrices, orbits of this group could be relatively easy to compute if  we have a convenient set of generators.  
For these codes, one may use an additional group $E\cong 2^{e}$ of sign changes at a space of codewords of $RM(1,d)$.
One can easily describe the action of $E$.  We therefore look for a union of orbits of $E(H\cap P)$.

We consider the action of $E(H\cap P)$ on the set of all $2^d$-tuples 
$\mathcal {A}$ of the form $(\pm 1, \pm 1, \dots , \pm 1)$ with respect to the standard sultry frame.  Also we have an action on 
$\mathcal {A}_0:=\mathcal A \cap \bw d$, which has cardinality $2^{1+d+{d \choose 2}}$.  

\begin{proc}\labttr{computer} 
Let $\mathcal O_i, i=1,2,\dots $ be the orbits of $E(H\cap P)$.  
An easy computer program can list these explicitly and compute inner products involving two orbits.  
  If only three different inner products occurs for some union $\mathcal O_j \cup \mathcal O_{j'} \cup \dots $ of two (or more!), this union is a tricosine spherical code.  Unlike in \reftt{unidefsc}, there is generally no reason to expect a transitive group of isometries.  
\end{proc}

\begin{rem} \labttr{app2} (i) 
A search for other codes could be done using other subgroups of $N$.  Since the $\bw d$ lattices contain vectors of shape 
$(1^X0^{\Omega \setminus X})$, for codewords $X \in RM(2,d)$, variations of $\mathcal A$ and $\mathcal A_0$ may be tried.

(ii) 
One can check whether the codes created this way are 
association schemes by straightforward  
accounting of inner products \reftt{computer}.  
\end{rem}

\section{The Optimism Code 
and a nonlinear $(16, 256, 6)$ binary code }

We shall define the {\it Optimism Code} or {\it Opticode}, a 4-cosine spherical code in dimension 16 with 256 unit vectors.  
A byproduct is that we deduce the existence of a 
nonlinear binary code with parameters $(16, 256, 6)$ and determine its automorphism group.  We furthermore deduce existence of a  64-point subcode with cosines $\{0, \pm \fourth\}$, which gives another existence proof of $\mathcal {NSC}_{16,64}$.

There is a famous nonlinear binary code with parameters $(16, 256, 6)$, the Nordstrom-Robinson code.  Existence of such a code has been given in several ways (see \cite{nr}, \cite{goe}).  There are references (e.g., in \cite{fst,ms}) to a uniqueness proof by S. L. Snover \cite{sno}, but the proof seems to be unpublished.   Our cocycle-style existence proof is probably new.

\subsection{Near-derivations for $AGL(4,2)$ on $RM(2,4)$}

The rest of this section uses the general discussion of $\bw 4$ and subgroups of its standard frame group, starting with \reftt{14setup}, but does not use the cohomology studies for $GL(3,2)$.  Instead, we study a much easier situation, that of degree 1 cohomology of $\FF_2[GL(4,2)]$ on its 6-dimensional module.  In fact, we really need only the more primitive concept of 1-cocycle.

\begin{nota}\labttr{neardergl42}
Let $M$ be the 6-dimensional irreducible module for $\FF_2[GL(4,2)]$ which occurs in the tensor square of the standard 4-dimensional module.  Then $H^1(GL(4,2),M)$ is 1-dimensional \cite{poll}.  If $X$ is the 8-dimensional permutation module for $GL(4,2)\cong Alt_8$, $X$ is uniserial with Loewey factors $\FF_2, M, \FF_2$ and every derivation on $M$ is inherited from  a derivation on $X$.  
\end{nota} 

\begin{nota}\labttr{ker} 
We use the notation of \reftt{neardergl42} and identify $D/(D\cap R)$ as a subquotient of $X$.  
Let $f$ be the near-derivation $P\rightarrow D$ whose associated derivation $\bar f$ 
on $D/(D\cap R)$ is identified with the derivation inherited from $X$ whose kernel $K$ is an  $2^4{:}Alt_7$ subgroup of $2^4{:}Alt_8$ (in more concrete language, we suppose that the permutation module for $Alt_8$  has basis $e_1, \dots , e_8$; then  $\bar f$ is identified with the map which sends permutation $g$ to $e_1-e_{1g}$ modulo $\FF_2(e_1+\dots + e_8)$).  
Then the set of nontrivial cosets of $D\cap R$ contained in $Im(f)(D\cap R)$ forms an orbit of length 7 for the action of $K=Ker(\bar f)$ \reftt{propderiv}.  
\end{nota}  

We need to check that $Im(f)(D\cap R)$ has weights 0, 6, 8, 10, 16 only (i.e., 4 and 12 do not occur).

\begin{lem}\labtt{ker4}   For the natural quadratic form on $RM(2,4)$, the radical is $RM(1,4)$.  The action of $AGL(4,2)$ on $RM(2,4)/RM(1,4)$ induces the associated $\Omega^+(6,2)$ and has kernel the translation subgroup.  
\end{lem} 
\pf
The first part follows from the well-known annihilation results for Reed-Muller codes.  The rest follows for example from group orders.  
\eop

\begin{coro} \labtt{ker5} The weights in $Im(f)(D\cap R)$ are just 0, 6, 8, 10 and  16.
\end{coro}   
\pf
Let $K_0$ be the stabilizer 
in $K$ of a coset of $D\cap R$ in the orbit of \reftt{ker}.  Then $K_0 \cong Alt_6$.  Such a coset in $D/(D\cap R)$ may be interpreted as a nonsingular vector in the sense of the natural nondegenerate  quadratic form on  $D/(D\cap R)$ (this is clear since the stabilizer of a singular vector in this orthogonal group is solvable).  
\eop

\begin{nota} \labttr{osc} 
We now let $\mathcal {OG}$ be the subgroup between $D\cap R$ and $N=DP$ \reftt{14setup} which corresponds to the near-derivation $f$ as in \reftt{ker}.  It has shape $2^{1+8}_+GL(4,2)$.  
Define 
$\mathcal {OC}:=x_0{\mathcal {OG}}$.  The stabilizer of $x_0$ in $\mathcal {OG}$ is just $\mathcal{OG} \cap P\cong 2^4{:}Alt_7$.  
This code clearly has cardinality 256 and the minus signs occur with multiplicities equal to  the weights of \reftt{ker5}.  

We call $\mathcal {OG}$ the {\it Optimism Group} and 
$\mathcal {OC}$ the {\it Optimism Code}.  For short, we say {\it Optigroup} and {\it Opticode}.  

The binary code $\mathcal {BC}_{16,256,6}$ is defined to be 
the set of 256 binary vectors corresponding to the elements of $\mathcal {OC}$ as follows: if $a=(a_i)\in \FF_2^{16}$ corresponds to $y=(y_i)\in \mathcal {OC}$, then $a_i= 0, 1$ according to whether $y_i=1, -1$, respectively.  
\end{nota}

\begin{nota}\labttr{diff}  For an orbit $X$, define $\DD (X):=\{ x-x'| x, x'\in X, x\ne \pm x'\}$.  
For a subset $A$ of $\Omega$, let $v_A:=\sum_{i \in A}v_i$ (so that $x_0=v_\Omega$).  
\end{nota}

\begin{lem}\labtt{ip1} 
(i) For any $B\subseteq \Omega$, 
$(x_0,x_0-2v_B)=16-2|B|$.

(ii) We have $(X,X)\subseteq \{0,  \pm 16\}$.

(iii) If $X, Y$ are different orbits, $x\in X, y \in Y$, then $(x,y)=\pm 4$.  
\end{lem}
\pf 
(i) 
We compute that $(x_0, v_B)=|B|$ and so $(x_0,x_0-2v_B)=16-2|B|$.   

(ii) This is clear since $x_0-2v_B$ is in the orbit of $x_0$ if and only if $B$ is a hyperplane or 0 or $\Omega$. 

(iii) We may assume by transitivity of $\mathcal {OG}$  that $x=x_0$. 
Then $y=x_0-2v_S$, where $S$ has defect 2, whence weight 6 or 10.  
\eop

\begin{lem} Let $M$ be the lattice spanned by all the vectors $\fourth \DD (X)$.  Then  $\fourth \DD (X)$ is the set of roots in $M$ and $M$ is a root lattice of type $D_{16}$.  
\end{lem}
\pf  By \reftt{ip1}(ii),  $\fourth \DD (X)$ is a root system with $480$ roots.  In dimension 16, the only candidates are $D_{16}$ and $E_8E_8$.  Since $M$ admits the irreducible action of $2^{1+8}_+Alt_7$, $E_8E_8$ is impossible (because, for example, this perfect group does not have an 8-dimensional representation with the center acting as $-1$).  
\eop

\begin{nota}\labttr{mx}  Suppose that $X$ is an orbit of $E$.  
Let $M'=M'(X)$ be the sublattice of $\dual M$ spanned by unit vectors in the dual, which are just the $\fourth x, x\in X$.  
Then $M'=M'(X)$ has isometry group $2\wr Sym_{16}$.  
\end{nota} 

\begin{lem} \labtt{uniserial} 
Suppose that $X=x_0E$.  The action of $EK\cong 2^{1+8}{:}Alt_7$ on $M'/2M'$ is uniserial with factors of dimensions 1, 4, 6, 4, 1. 
\end{lem} 
\pf  The action is that of $K\cong 2^4{:}Alt_7$, which is contained in $P$, the natural $AGL(4,2)$-subgroup of $Sym_{16}$.  The factors are inherited from the action of a natural $GL(4,2)$-subgroup and uniseriality follows from  commutation by the normal subgroup of order 16, which is in $EK$.  

The last statement follows from the structure of the dual of the $D_{16}$-lattice and \reftt{ip1}(i).  
\eop

\begin{lem} \labtt{xrm}  
Let $X$ be an orbit, $M'=M'(X)$.   

(i) 
Let $y \in Y$, an orbit different from $X$.  Then 
$M'_y:=\{v\in M' | (v,y)\in 2\ZZ \}$ has index 2 in $M'$.  
Suppose $X=x_0E$.  
The span of the image of all $Y\ne X$ in $Hom(M',\ZZ_2)$ with respect to the double  basis $\fourth  X$ corresponds to $RM(2,4)$.  

(ii) Furthermore, when $X=x_0E$, 
$\cap_{y \in Y, Y\ne X} M'_y$ is 
a sublattice between $2M'$ and $M'$ which corresponds to $RM(1,4)$ with respect to the basis $\fourth X$ modulo $2M'$.  

(iii) Let $T_X$ be the normal subgroup of 
$Stab_{Aut(\mathcal {OC})}(X)$ 
corresponding to the group of sign changes \reftt{mx}.    
Then, for orbits $Y\ne X$,  $T_X\cap T_Y=\la -1 \ra$.  
\end{lem} 
\pf 
(i) and (ii) are  clear from \reftt{ip1} and the well-known annihilation result 
$RM(2,4)^\perp=RM(1,4)$.  

(iii) We may assume that 
$X=x_0E \ne Y$.  
Let $t\in T_X\cap T_Y$.    Given $x\in X$, there is a scalar, $c$, so that $x^t=cx$.  Similarly, given $y \in Y$, there exists a scalar, $d$, so that $y^t=dy$.  By \reftt{ip1}, $0 \ne (x,y)=(x^t,y^t)=(cx,dy)=cd(x,y)$, whence $cd=1$.  Since $c, d \in \{ \pm 1\}$, $c=d$.   Therefore $t$ acts as the scalar $c$.  
\eop

\begin{coro}\labtt{normalstabx}  
Let $X:=x_0E$.  Then 
$Stab_{Aut(\mathcal {OC})}(X)=EK \le \brw 4$. 
\end{coro} 
\pf 
By \reftt{xrm}(iii), $T_Y/\la -1\ra$ embeds in 
$Stab_{Aut(\mathcal {OC})}(X)/T_X$ as a subgroup stabilized by the action 
$Stab_{Aut(\mathcal {OC})}(X) \cap Stab_{Aut(\mathcal {OC})}(Y)$, which has a quotient isomorphic to $Sym_6$.   Since 
$Stab_{Aut(\mathcal {OC})}(X)/T_X\cong 2^4{:}Alt_7$, it follows that 
  $T_Y/\la -1\ra$ has order dividing $2^4$.  Therefore, 
  $Stab_{Aut(\mathcal {OC})}(X)$ has order dividing $2^9|Alt_7|$.  Since $Stab_{Aut(\mathcal {OC})}(X)$ contains $EK\cong 2^{1+8}Alt_7$, we have equality.  
\eop

\begin{thm}\labtt{autopt}  
The isometry group of the optimism code is just 
the optimism group $\mathcal {OG}$, 
of shape  $2^{1+8}_+GL(4,2)$.    
\end{thm}
\pf  
Use \reftt{normalstabx} and the fact that $\mathcal {OG}$ is transitive on the orbits of $E$. 
\eop  
\begin{prop}\labtt{autbc} The isometry group of $\mathcal{BC}_{16,256}$ is $2^4{:}Alt_7$.  
\end{prop}
\pf
The isometry group of this binary code embeds 
by coordinatewise action in 
$Aut(\mathcal  {OSC}_{16,256})$ as the subgroup 
$Aut(\mathcal  {OSC}_{16,256})\cap P$ 
of the isometry group of $\mathcal  {OC}_{16,256}$.  This is the subgroup stabilizing $x_0$.  
\eop 

\begin{rem} \labttr{historyautbc}  For earlier determinations of the automorphism group of the Nordstrom-Robinson code, see \cite{ber}, \cite{goe}, \cite{sno}. 
\end{rem}

\subsection{$\mathcal {NSC}_{16,64}$ as subcode of the Optimism Code} 

. 

\begin{nota}\labttr{nagl42} We let $X$ be the subgroup between $E$ and $\mathcal {OG}$ which corresponds to the near-derivation $f$ restricted to  $P_0$, the stabilizer in $P$ of $0$.  So, $P_0\cong GL(4,2)$ and $X \cong 2^5.GL(4,2)$. 

Let $E_{01}:=C_E(W_{01})$ be the rank 3 subgroup of $E$ which is trivial on $W_{01}$ (see \reftt{14setup}).  Then $N_X(E_{01})$ has the form $2^5.2^3.GL(3,2)$ and 
$N_P(E_{01})=P\cap N_X(E_{01})\cong GL(3,2)$, which acts indecomposably on 
$E/Z(R)$.  Define $H:=N_X(E_{01})\cap C(W_{01})$, a subgroup of index 4 in $N_X(E_{01})$ which satisfies $H\cap R=E_{01}$.  Then $H$ has shape $2^3.2^3.GL(3,2)$, $H/E_{01}\cong AGL(3,2)$.  
Note that $p_{01}$ fixes $x_0$ and normalizes $H$; in fact, $[H,p_{01}]=E_{01}$.    
\end{nota}

\begin{lem}\labtt{nonsplit} The group $H$ of \reftt{nagl42} has shape $4^3{:}GL(3,2)$.
\end{lem} 
\pf
The action of an element of order 7 in $H$ on $O_2(H)$ forces $O_2(H)$ to be abelian, since the two composition factors therein are isomorphic (by commutation with $p_{01}$, for example).  
We assume that  $O_2(H)$ is elementary abelian, then derive a contradiction.

Then $O_2(H)$ is completely reducible as a module for $H/O_2(H)\cong GL(3,2)$ (it is easy to prove that $Ext^1(Y,Y)=0$ for a 3-dimensional irreducible $Y$ with a minimal resolution \cite{benson}).  
 There are subgroups isomorphic to 
$GL(3,2)$ in $X\cap P\cong 2^4{:}Alt_7$, acting indecomposably on $\Omega$, fixing $\{0, 1\}$.  Therefore, $H$ splits over $H\cap R=E_{01}$.  By Gasch\"utz's theorem, $X$ splits over $O_2(X)$.
This is a contradiction to \reftt{summary1}, \reftt{nonsplitgl32} since the nonsplit extension $2^3{\cdot}GL(3,2)$ does not embed in $AGL(4,2)$.  

We conclude that  $O_2(H) \cong 4^3$.  
\eop

\begin{de}\labttr{seconddefbcghcode}  We define 
the spherical code $\mathcal S :=\fourth x_0 H$, where $H$ is as in \reftt{nagl42}.  The cosines are just $\{0, \pm \fourth \}$.  Its cardinality is $|H:H\cap P|=64$.  In fact, $O_2(H)$ acts regularly on $\mathcal S$.  
\end{de} 

\begin{rem}\labttr{recover} (i) 
We may view $H$ as the subgroup between $E$ and $N_X(E_{01})$ which corresponds to the near-derivation $f$ restricted to the subgroup of $P_0$ which stabilizes the points $\o_0$ and $\o_1$, equivalently, which normalizes $E_{01}$.
Therefore, $\mathcal S$ is identified with 
$\mathcal{NSC}_{16,64}$, which was defined as $x_0H^*$, where $H^*$ is the subgroup  in \reftt{x0} (reason: the cocycle $f$ we used in \reftt{ker} could have been used to define a suitable group $H^*$ as in \reftt{x0} since its values have the right weights).  

(ii) This new realization of $\mathcal{NSC}_{16,64}$ has the advantage of exhibiting a larger group of isometries, $H\la p_{01} \ra$, 
than $H^*\la p_{01}\ra$.  Upon projection to 14-space $W$, we get a code like $\mathcal {BCGM}$.

(iii) For another discussion of our extensions, see \cite{grdemp}.  
\end{rem}

\subsection{Concluding Remarks}

\begin{rem} \labttr{256code}{\it Alternate constructions of the Optimism Code and a $(16, 256, 6)$ nonlinear binary code.}    We take 
our spherical code 
$\mathcal {NSC}_{16, 64}$ and the group $D\cap R\cong 2^5$.  The new spherical code $\mathcal {OSC}_{16, 256} :=(\mathcal {NSC}_{16, 64})(D\cap R)$ has 256 vectors and cosine set 
$\{0,\pm \fourth, -1\}$.   
The binary code $\mathcal {BC}_{16,256,6}$ is a set of 256 binary vectors corresponding to the elements of $\mathcal {NSC}_{16, 256}$ as follows: $a=(a_i)\in \FF_2^{16}$ corresponds to $y=(y_i)\in \mathcal {NSC}_{16, 256}$ by the rule $a_i= 0, 1$ according to whether $y_i=1, -1$, respectively.  
Finally, one may start with a Nordstrom-Robinson type binary code and reverse the previous procedure to define a spherical code.  
\end{rem}  

\begin{rem}\labttr{sce} {\it Spherical codes and energy.}  
It is clear that one can make many spherical codes in $\RR^n$ by taking orbits of the all-1 vector by subgroups of the degree $n$ monomial group $Mon(n, \{\pm 1\})$.  
One can get larger spherical codes as orbits by overgroups of such monomial groups, e.g. the optimism group is contained in a natural 
$2^{1+8}.Alt_9$ subgroup of $\brw 4$. 
There are many  candidates to try.  
It is not clear which are likely to be associated to universally optimal situations.  
Known examples involve exceptional objects as well as series (see the table on page 2 of \cite{ck}).    
\end{rem}

\section{Appendix: Background on 1-cocycles and derivations} 

\begin{de}\labttr{rightderiv} A {\it right 1-cocycle} or {\it right derivation} from the group $X$ to the additive right $X$-module $A$ is a function $f:X\rightarrow A$ so that 
$f(xy)=f(x)+f(y)^{x^{-1}}$ for all $x, y \in X$.  
A {\it 1-coboundary} or  {\it inner derivation} is such a function of the form $f(x)=a-ax$, for a fixed $a \in A$.  A noninner derivation is sometimes called an {\it outer derivation}.  

In case $A$ is a multiplicative group, the derivation condition reads 
$f(xy)=f(x)f(y)^{x^{-1}}$. The inner derivation condition reads $f(x)=aa^{-x}$.   
\end{de}  

\begin{prop}\labtt{interp1coc}  Let the group $H$ be a semidirect product of normal abelian subgroup $A$ by a complement $X$.  The complements correspond to the 1-cocycles from $X$ to $A$:  if $f$ is a 1-cocycle, the complement associated to it is 
$\{f(x)x|x \in H\}$.   
Two complements are conjugate by $H$ (equivalently, by $A$) if and only if their corresponding 1-cocycles are cohomologous (i.e., their difference is a 1-coboundary).
\end{prop}
\pf  Classic.  See for example \cite{gruenberg}, \cite{hup}.  \eop

\begin{de}\labttr{kerder} The {\it kernel} of a derivation as in \reftt{rightderiv} is $Ker(f):=\{g\in X|f(g)=0\}$ (in the additive case) and
$Ker(f):=\{g \in X| f(g)=1\}$ (in the multiplicative case).  It is a subgroup, though typically not normal.
\end{de}

\begin{lem} \labtt{propderiv}
Let $f:X\rightarrow A$ be as in \reftt{rightderiv}.  Let $K$ be the kernel of the derivation $f$.  Then

(i) $f$ is constant on left cosets of $K$ in $X$; also, if $x, y \in K$, then $f(x)=f(y)$ if and only if $xK=yK$; 

(ii) If $x\in K$, $f(xy)=f(y)^{x^{-1}}$; consequently, the values of $f$ on the right coset $Ky$ of $K$ in $X$ form a $K$-orbit in $A$.

(iii) The values of $f$ on the double coset $KxK$, for $x\in X$, form a $K$-orbit on $A$, the orbit containing $f(x)$.  

(iv)  Suppose $X$ acts doubly transitively on the cosets of $K$.  Then the set of values of $f$ is the disjoint union of $0\in A$ with the $K$-orbit of values taken by $f$ on the nontrivial double coset of $K$ in $X$.  

\end{lem}  
\pf 
Easy work with the definition of derivation.  For (i), set $a=xy, b=x$.  Then $f(a)=f(b)f(b^{-1}a)^{b^{-1}}$.  Consider the condition $aK=bK$.  
\eop

\begin{prop}\labtt{1coh168} 
(i) $dim(H^1(GL(n,2),\FF_2^n))=0$ if $n \ne 3$; 

(ii)  
$dim(H^1(GL(3,2),\FF_2^3))=1$. 
\end{prop}  
\pf 
(i) and (ii) may be found in \cite{higman}.

The result (ii) 
is well-known and follows trivially from modular representation theory, specifically the structure of projective indecomposable modules for $\FF_2[GL(3,2)]$.  For a proof with resolutions, see \cite{benson}.   For an elementary proof using the interpretation of complements modulo conjugacy, see \cite{hup} \cite{g12}.  For another, see the proof of \reftt{lowdegcoh168}.  
\eop 

\begin{lem}\labtt{anatomy}
Let $G\cong GL(3,2)$ and $M$ a 3-dimensional irreducible $\FF_2G$-module.  

(i) If $f$ is a nonzero inner derivation, then $Im(f)$ is a 7-subset of $M$ containing 0.  If $a\in M\setminus Im(f)$, then $f$ is the inner derivation $x\mapsto a(x-1)$.  

(ii) If $f$ is a noninner derivation, $Im(f)=M$.  Also, $Ker(f)\cong 7{:}3$, the Frobenius group of order 21.  
\end{lem}
\pf
(i)  Since $f$ is inner, there is $a \in M$ so that $f(x)=a(1-x)$.  The kernel of $f$ is the index 7 stabilizer of $a$ ($a\ne 0$ since $f\ne 0$).  Obviously there is no solution to $f(x)=a$.  Now use \reftt{propderiv}(i).  

(ii)  
There is an indecomposable module, $L$, so that $M\le L$ and $L/M\cong \FF_2$.  There is $a\in L$ so that $f(x)=a(1-x)$ for all $x \in G$.  Since $f$ is noninner, $a \in L\setminus M$.  The stabilizer $K$ of $a$ in $G$ (equivalently, the kernel of the derivation $f$) must have even index, or by Maschke's theorem, $L$ would be decomposable.  

The kernel of $f$ is a proper subgroup with index at most $|M|=8$.  
The index is therefore 7 or 8.    
The last paragraph proves the index is not 7,  and so we are done by 
\reftt{propderiv}(i).
\eop

\begin{prop}\labtt{lowdegcoh168} 
Let $G\cong GL(3,2)$ and let $M$ be a 6-dimensional indecomposable module for $\FF_2G$ with composition factors of dimension 3 and which are duals of each other.  
Let $S$ be the socle of $M$.  
Write 
$0\rightarrow S \rightarrow M \rightarrow M/S \rightarrow 0$.

Then 

(i) $H^1(G,M)$ has dimension 1 and is the image of $H^1(G, S)$ under the natural map coming from $S\rightarrow M$;

(ii) Let  $f:G\rightarrow M$ be a 1-cocycle.  Either (a) the values of $f$ are contained in $S$; or (b), the values of $f$ 
are not contained in $S$  and the kernel of $f$ is contained in the index 7 subgroup of $G$ which stabilizes a nonzero vector of $M/S$.  
\end{prop}  
\pf 
(i) 
Let $P$ be the permutation module for $G$ on the cosets of a subgroup $H$ of index 7.  There are two conjugacy classes of such $H$ and for one of them, $P$ is isomorphic to the direct sum of $M$ and the trivial module.  We have $H^1(G,P)\cong H^1(H,\FF_2)$ by the Eckmann-Shapiro lemma.  The right object  has dimension 1 since $H\cong Sym_4$ and $H^1(H,\FF_2)\cong Hom(H, \FF_2)$.  By additivity, $H^1(G,P)\cong H^1(G,M)\oplus H^1(G,\FF_2)$ and the last summand is zero since it is isomorphic to $Hom(G,\FF_2)$.  This proves that $H^1(G,M)$ has dimension 1.  

From the long exact sequence, we get 
$0\rightarrow H^1(G,S)\rightarrow H^1(G,M)\rightarrow H^1(G,M/S)\rightarrow \dots$.  By dimensions, using the preceding paragraph and \reftt{1coh168}(ii), we get the final statement.  

Actually, we can prove 
\reftt{1coh168}(ii) directly.  This long exact sequence proves that at least one of 
$H^1(G,S), H^1(G,M/S)$ is nonzero.  It follows that both are nonzero since the modules $S$ and $M/S$ are related by an outer automorphism of $G$.  The long exact sequence then proves that one has dimension 1, so both have dimension 1.  

(ii)  We may assume that  the values of $f$ do not lie in $S$.  By (i), the composition of $f$ with the quotient modulo $S$ is cohomologous to 0, i.e., there is $a\in M$ so that $f(x) - (a-ax) \in S$ for all $x \in G$.  Such an $a$ is not in $S$.  The kernel of $f$ is therefore contained in the stabilizer in $G$ of the nontrivial vector $a+S \in M/S$. 
\eop

\begin{lem}\labtt{involext} 
Let $M$ be an irreducible 3-dimensional module for $\FF_2G$, where $G\cong GL(3,2)$ and let $H$ be an extension, so that $M$ is normal in $H$ and $H/M\cong G$.

Then 

(i) $Aut(H)$ is an extension of $Inn(H)$ by $\la u \ra$, where the involution $u$ acts trivially on $M$ and on $H/M$, so induces a noninner  derivation from $H/M$ to $M$;

(ii) $O_2(Aut(H))$ acts transitively on the two $H$-classes of involutions in $H\setminus M$; moreover, if $t_1$ and $t_2$ are involutions so that $Mt_1=Mt_2$, there exists $g \in O_2(Aut(H))$ so that $g$ takes $t_1$ to $t_2$ (more precisely, if $t\in G$ is in $Mt_1$, there exists a derivation $f: G\rightarrow M$ so that $f(t)=t_1t_2$).  

(iii) If $H$ is a split extension and $t\in H \setminus M$ is an involution, there exists a complement to $M$ in $H$ which contains $t$.  
\end{lem}
\pf
(i) Since 
$Aut(G)$ interchanges the two irreducible 3-dimensional $G$-modules, $Aut(H)$ induces only $Inn(G)$ on $G\cong H/M$.  Let $R$ be $C_{Aut(H)}(H/M)$.  Then $R$ acts trivially on $M$, by absolute irreducibility of $M$.  Then $R$ is identified with the 1-cocycles from $G$ to $M$, which forms a 4-dimensional space, by \reftt{1coh168}.  

(ii) If $H$ is a split extension, $H\setminus M$ contains involutions, and if it is nonsplit, the same is true, by \cite{gnw}, for example.  

By looking at Jordan canonical forms, it is clear that $H$ has two conjugacy classes of involutions.  The kernel of an outer derivation is a Frobenius group of order 21.  Therefore, an involution in $H\setminus M$ has, in its action on $O_2(Aut(H))$, Jordan canonical form which is a sum of two indecomposable $2\times 2$ blocks.  Therefore, the two $H$-classes of such involutions fuse in $Aut(H)$.  
Since the space of derivations is a free module for $\la t \ra$, an involution in $Mt_1$, every fixed point is the image of a derivation.  Take the fixed point $t_1t_2$ to prove the final statement.

(iii) This follows from (ii).   
\eop 

\begin{lem}\labtt{tensorgl32}
$\threed \otimes \threed$ is uniserial and has composition factors $\threedp, \threed, \threedp$; 
$\threedp \otimes \threedp$ is uniserial and has composition factors $\threed, \threedp, \threed$; 
$\threed \otimes \threedp \cong \oned \oplus \eightd$.
\end{lem}
\pf 
Well known, and easy to prove with Brauer characters. 
\eop 

\begin{lem}\labtt{nonsplitgl32}
Suppose that $1\rightarrow A \rightarrow 
E \rightarrow G \rightarrow 1$ is an extension of $G\cong GL(3,2)$ by its standard module $A \cong \FF_2^3$.  Let $A_1$ be a maximal subgroup of $A$ and let $C_1:=C_E(A)$.  Then either $C_1$ splits over $A$ or there exists $B_1 \le C_1, B_1 \cong 4^2$ and the elements of $C_1\setminus B_1$ invert $B_1$ under conjugation.
\end{lem}
\pf
Let $N_1:=N_E(A_1)$.  Then $N_1/C_1\cong GL(2,2)\cong Sym_3$.  Clearly, $C_1$ is self-centralizing in $E$ and is homocyclic abelian.  
Define $B_1:=[C_1, N_1]=[C_1,O_3(N_1)]$.  Either $B_1$ is elementary abelian or   
$B_1 \cong 4^2$ and the elements of $C_1\setminus B_1$ invert $B_1$ under conjugation.  In either case, $C_1\setminus A$ contains involutions and so every coset $xA\ne A$ which satisfies $x^2\in A$ contains involutions.

Suppose that $B_1$ is elementary abelian.  Let $t$ be an involution in $N_E(B_1)\setminus O_2(N_E(B_1))$.  Then $t$ inverts $h$, an element of order 3 in $N_E(B_1)$, by the Baer-Suzuki theorem.  The module 2-dimensional faithful module for $\FF_2[Sym_3]$ is projective and injective, so $B_1$ has a splitting $A_1 \times A_2$ as modules for $\la t, h \ra$.  Then $A_2\la t\ra\cong Dih_8$ meets $A$ trivially.  By Gash\"utz's theorem, $E$ splits over $A$. 

If $B_1$ is not elementary abelian, the second alternative.  \eop

\end{document}